\documentclass[12pt]{amsart}
\usepackage{latexsym,amssymb,amsmath,amsthm,amscd,enumerate,graphicx,esint,longtable}

\makeatletter
\@namedef{subjclassname@2010}{%
  \textup{2010} Mathematics Subject Classification}
\makeatother

\newtheorem{thm}{Theorem}[section]
\newtheorem{cor}[thm]{Corollary}
\newtheorem{lem}[thm]{Lemma}

\newtheorem{prob}[thm]{Problem}

\theoremstyle{definition}

\newtheorem*{xrem}{Remark}

\numberwithin{equation}{section}

\newcommand{\cA}{{\mathcal A}}
\newcommand{\cB}{{\mathcal B}}

\newcommand{\cQ}{{\mathcal Q}}
\newcommand{\cR}{{\mathcal R}}

\newcommand{\R}{{\mathbb R}}

\newcommand{\Z}{{\mathbb Z}}

\newcommand{\kM}{{\mathfrak M}}

\newcommand{\kT}{{\mathfrak T}}

\def\al{\alpha}
\def\bt{\beta}

\def\dl{\delta}

\def\Sg{\Sigma}
\def\om{\omega}

\def\0{\emptyset}
\def\1{\textbf{\rm 1}}
\def\6{\partial}
\def\8{\infty}

\def\lt{\left}
\def\rt{\right}

\def\wt{\widetilde}

\DeclareMathOperator*{\essinf}{\,\text{ess\,inf}\,}

\begin{document}

\title
[The Fefferman-Stein type inequalities]
{The Fefferman-Stein type inequalities for strong and directional maximal operators in the plane}

\author[H.~Saito]{Hiroki~Saito}
\address{
Academic support center, 
Kogakuin University, 
2665--1, Nakanomachi, Hachioji-shi Tokyo, 192--0015, Japan
}
\email{j1107703@gmail.com}
\author[H.~Tanaka]{Hitoshi~Tanaka}
\address{
Research and Support Center on Higher Education for the hearing and Visually Impaired, 
National University Corporation Tsukuba University of Technology,
Kasuga 4-12-7, Tsukuba City, Ibaraki, 305-8521 Japan
}
\email{htanaka@k.tsukuba-tech.ac.jp}

\thanks{
The first author is supported by 
Grant-in-Aid for Young Scientists (B) (15K17551), 
the Japan Society for the Promotion of Science. 
The second author is supported by 
Grant-in-Aid for Scientific Research (C) (15K04918), 
the Japan Society for the Promotion of Science. 
}

\subjclass[2010]{Primary 42B25; Secondary 42B35.}

\keywords{
directional maximal operator;
Fefferman-Stein type inequality;
Hardy-Littlewood maximal operator;
strong maximal operator.
}

\date{}

\begin{abstract}
The Fefferman-Stein type inequalities 
for strong maximal operator and directional maximal operator 
are verified with composition of the Hardy-Littlewood maximal operator in the plane.
\end{abstract}

\maketitle

\section{Introduction}\label{sec1}
The purpose of this paper is to develop a theory of weights for 
strong maximal operator and directional maximal operator in the plane. 
We first fix some notations. 
By weights we will always mean non-negative and locally integrable functions on $\R^n$. 
Given a measurable set $E$ and a weight $w$, 
$w(E)=\int_{E}w(x)\,dx$, 
$|E|$ denotes the Lebesgue measure of $E$ and 
$\1_{E}$ denotes the characteristic function of $E$. 
Let $0<p\le\8$ and $w$ be a weight. 
We define the weighted Lebesgue space 
$L^p(\R^n,w)$ 
to be a Banach space equipped with the norm (or quasi norm) 
\[
\|f\|_{L^p(\R^n,w)}
=
\lt(\int_{\R^n}|f(x)|^pw(x)\,dx\rt)^{1/p}.
\]
For a locally integrable function $f$ on $\R^n$, 
we define the Hardy-Littlewood maximal operator $\kM_{\cQ}$ by
\[
\kM_{\cQ}f(x)
=
\sup_{Q\in\cQ}
\1_{Q}(x)\fint_{Q}|f(y)|\,dy,
\]
where $\cQ$ is the set of all cubes in $\R^n$ 
(with sides not necessarily parallel to the axes) 
and the barred integral $\fint_{Q}f(y)\,dy$ 
stands for the usual integral average of $f$ over $Q$. 
For a locally integrable function $f$ on $\R^n$, 
we define the strong maximal operator $\kM_{\cR}$ by
\[
\kM_{\cR}f(x)
=
\sup_{R\in\cR}
\1_{R}(x)\fint_{R}|f(y)|\,dy,
\]
where $\cR$ is the set of all rectangles in $\R^n$ 
with sides parallel to the coordinate axes. 

Let $\kT:\,L^p(\R^n)\to L^p(\R^n)$, $p>1$, 
be a sublinear operator. 
It is a fundamental problem of the weight theory that 
to determine some maximal operator $\kM_{\kT}$
capturing certain geometric characteristics of $\kT$ such that 
\begin{equation}\label{1.1}
\int_{\R^n}|\kT f(x)|^pw(x)\,dx
\le C
\int_{\R^n}|f(x)|^p\kM_{\kT}w(x)\,dx
\end{equation}
holds for arbitrary weight $w$. 
It is well known that 
\[
\int_{\R^n}\kM_{\cQ}f(x)^pw(x)\,dx
\le C
\int_{\R^n}|f(x)|^p\kM_{\cQ}w(x)\,dx
\]
holds for arbitrary weight $w$ and $p>1$,
and further that
\begin{equation}\label{1.2}
\sup_{t>0}
t\,w(\{x\in\R^n:\,\kM_{\cQ}f(x)>t\})
\le C
\int_{\R^n}|f(x)|\,\kM_{\cQ}w(x)\,dx.
\end{equation}
These are called the Fefferman-Stein inequality 
and are toy models of the problem \eqref{1.1} 
(see \cite{FS}).

There is a problem in the book \cite[p472]{GR}:

\begin{prob}\label{prob1.1}
Does the analogue of the Fefferman-Stein inequality 
hold for the strong maximal operator, i.e. 
\begin{equation}\label{1.3}
\int_{\R^n}\kM_{\cR}f(x)^pw(x)\,dx
\le C
\int_{\R^n}|f(x)|^p\kM_{\cR}w(x)\,dx,
\quad p>1,
\end{equation}
for arbitrary $w(x)\ge 0$?
\end{prob}

\noindent
Concerning Problem \eqref{prob1.1} 
it is known that, see \cite{Lin} (also \cite{Pe1,Pe2}), 
\eqref{1.3} holds for all $p>1$ 
if $w\in A_{\8}^{*}$. 

We say that $w$ belongs to 
the class $A_p^{*}$ whenever
\begin{align*}
[w]_{A_p^{*}}
&=
\sup_{R\in\cR}
\fint_{R}w(x)\,dx\lt(\fint_{R}w(x)^{-1/(p-1)}\,dx\rt)^{p-1}
<\8,
\quad 1<p<\8,
\\
[w]_{A_1^{*}}
&=
\sup_{R\in\cR}
\frac{\fint_{R}w(x)\,dx}{\essinf_{x\in R}w(x)}
<\8.
\end{align*}
It follows by H\"{o}lder's inequality that 
the $A_p^{*}$ classes are increasing, 
that is, 
for $1\le p\le q<\8$ we have 
$A_p^{*}\subset A_q^{*}$. 
Thus one defines 
\[
A_{\8}^{*}
=
\bigcup_{p>1}A_p^{*}.
\]

The endpoint behavior of $\kM_{\cR}$ close to $L^1$ is given by 
Mitsis \cite{Mi} (for $n=2$) and 
Luque and Parissis \cite{LP} (for $n>2$). 
That is, 
\begin{align*}
\lefteqn{
w(\{x\in\R^n:\,\kM_{\cR}f(x)>t\})
}\\ &\le C
\int_{\R^n}
\frac{|f(x)|}{t}
\lt(1+\lt(\log^{+}\frac{|f(x)|}{t}\rt)^{n-1}\rt)
\kM_{\cR}w(x)\,dx,
\quad t>0,
\end{align*}
holds for any $w\in A_{\8}^{*}$, 
where $\log^{+}t=\max(0,\log t)$. 

In this paper concerning Problem \eqref{prob1.1} 
we shall establish the following. 

\begin{thm}\label{thm1.2}
Let $w$ be any weight on $\R^2$ and 
set $W=\kM_{\cR}\kM_{\cQ}w$. 
Then 
\begin{align*}
\lefteqn{
w(\{x\in\R^2:\,\kM_{\cR}f(x)>t\})
}\\ &\le C
\int_{\R^2}
\frac{|f(x)|}{t}
\lt(1+\log^{+}\frac{|f(x)|}{t}\rt)
W(x)\,dx,\quad t>0,
\end{align*}
holds, where the constant $C>0$ does not depend on $w$ and $f$.
\end{thm}

By interpolation, we have the following corollary.

\begin{cor}\label{cor1.3}
Let $w$ be any weight on $\R^2$ and 
set $W=\kM_{\cR}\kM_{\cQ}w$. 
Then, for $p>1$, 
there exists a constant $C_p>0$ such that 
\[
\|\kM_{\cR}f\|_{L^p(\R^2,w)}
\le C_p
\|f\|_{L^p(\R^2,W)}
\]
holds for all $f\in L^p(\R^2,W)$.
\end{cor}

Let $\Sg$ be a set of unit vectors in $\R^2$, 
i.e., a subset of the unit circle $S^1$. 
Associated with $\Sg$, 
we define $\cB_{\Sg}$ to be the set of all rectangles in $\R^2$ 
whose longest side is parallel to some vector in $\Sg$. 
For a locally integrable function $f$ on $\R^2$, 
we also define the directional maximal operator $\kM_{\Sg}$ 
associated with $\Sg$ as 
\[
\kM_{\Sg}f(x)
=
\sup_{R\in\cB_{\Sg}}
\1_{R}(x)\fint_{R}|f(y)|\,dy.
\]
Many authors studied this operator, 
see \cite{ASV1,ASV2,Ka1,Ka2,ST,Str},
and Katz showed that $\kM_{\Sg}$ is bounded on $L^2(\R^2)$ 
with the operator norm $O(\log N)$ 
for any set $\Sg$ with cardinality $N$.

For fixed sufficiently large integer $N$, let 
\[
\Sg_{N}
=
\lt\{\lt(
\cos\frac{2\pi j}{N},
\sin\frac{2\pi j}{N}
\rt):\,
j=0,1,\ldots,N-1
\rt\}
\]
be the set of $N$ uniformly spread directions on the circle $S^1$.
In this paper we shall prove the following,
which is a weighted version of the classical result due to 
Str\"{o}mberg \cite{Str}.

\begin{thm}\label{thm1.4}
Let $N>10$ and $w$ be any weight on $\R^2$. 
Set $W=\kM_{\Sg_{N}}\kM_{\cQ}w$. 
Then 
\begin{equation}\label{1.4}
\sup_{t>0}
t\,w(\{x\in\R^2:\,\kM_{\Sg_{N}}f(x)>t\})^{1/2}
\le C(\log N)^{1/2}
\|f\|_{L^2(\R^2,W)}
\end{equation}
holds for all $f\in L^2(\R^2,W)$, 
where the constant $C>0$ does not depend on $w$ and $f$.
\end{thm}

By interpolation, we have the following corollary.

\begin{cor}\label{cor1.5}
Let $N>10$ and $w$ be any weight on $\R^2$. 
Set $W=\kM_{\Sg_{N}}\kM_{\cQ}w$. 
Then, for $2<p<\8$, 
there exists a constant $C_p>0$ such that 
\[
\|\kM_{\Sg_{N}}f\|_{L^p(\R^2,w)}
\le C_p(\log N)^{1/p}
\|f\|_{L^p(\R^2,W)}
\]
holds for all $f\in L^p(\R^2,W)$.
\end{cor}

The letter $C$ will be used for the positive finite constants that may change from one occurrence to another. 
Constants with subscripts, such as $C_1$, $C_2$, do not change
in different occurrences.

\section{Proof of Theorem \ref{thm1.2}}\label{sec2}
In what follows we shall prove Theorem \ref{thm1.2}. 
Our proof relies upon the refinement of the arguments in \cite{Mi}. 
With a standard argument, we may assume that 
the basis $\cR$ is the set of all dyadic rectangles $R$ 
(cartesian products of dyadic intervals) 
having long side pointing in the $x_1$-direction. 
We denote by $P_i$, $i=1,2$, 
the projection on the $x_i$-axis. 
Fix $t>0$ and 
given the finite collection of dyadic rectangles 
$\{R_i\}_{i=1}^{M}\subset\cR$ 
such that 
\begin{equation}\label{2.1}
\fint_{R_i}|f(y)|\,dy>t,\quad
i=1,2,\ldots,M.
\end{equation}
It suffices to estimate 
$w\lt(\bigcup_{i=1}^{M}R_i\rt)$ 
(see the next section for details).

First relabel if necessary so that 
the $R_i$ are ordered so that 
their long sidelengths $|P_1(R_i)|$ decrease. 
We now give a selection procedure to find subcollection 
$\{\wt{R}_i\}_{i=1}^{N}
\subset
\{R_i\}_{i=1}^{M}$. 

Take $\wt{R}_1=R_1$ and 
let $\wt{R}_2$ be the first rectangle $R_j$ such that 
\[
|R_j\cap\wt{R}_1|<\frac13|R_j|.
\]
Suppose have now chosen the rectangles 
$\wt{R}_1,\wt{R}_2,\ldots,\wt{R}_{i-1}$. 
We select $\wt{R}_i$ to be the first rectangle $R_j$ 
occurring after $\wt{R}_{i-1}$ 
so that 
\[
\lt|\bigcup_{k=1}^{i-1}R_j\cap\wt{R}_k\rt|<\frac13|R_j|.
\]
Thus, we see that 
\begin{equation}\label{2.2}
\lt|\bigcup_{j=1}^{i-1}\wt{R}_i\cap\wt{R}_j\rt|<\frac13|\wt{R}_i|,
\quad i=2,3,\ldots,N.
\end{equation}

We claim that 
\begin{equation}\label{2.3}
\bigcup_{i=1}^{M}R_i
\subset
\lt\{
x\in\R^2:\,
\kM_{\cQ}[\1_{\bigcup_{i=1}^{N}\wt{R}_i}](x)
\ge\frac13
\rt\}.
\end{equation}
Indeed, 
choose any point $x$ inside a rectangle $R_j$ 
that is not one of the selected rectangles $\wt{R}_i$. 
Then, 
there exists a unique $K\le N$ such that 
\[
\lt|\bigcup_{i=1}^{K}R_j\cap\wt{R}_i\rt|
\ge\frac13|R_j|.
\]
Since, 
$|P_1(R_j)|\le|P_1(\wt{R}_i)|$ 
for $i=1,2,\ldots,K$, 
we have 
\[
P_1(R_j)\cap P_1(\wt{R}_i)
=
P_1(R_j)
\text{ when }
R_j\cap\wt{R}_i\ne\0,
\]
where we have used the dyadic structure;
\begin{equation}\label{2.4}
\text{
If both $I$ and $J$ are the dyadic interval then 
$I\cap J\in\{I,J,\0\}$.
}\end{equation}
Thus, 
\[
\bigcup_{i=1}^{K}
R_j\cap\wt{R}_i
=
\bigcup_{i=1}^{K}
P_1(R_j)
\times
(P_2(R_j)\cap P_2(\wt{R}_i))
=
P_1(R_j)
\times
\bigcup_{i=1}^{K}
P_2(R_j)\cap P_2(\wt{R}_i).
\]
Hence, 
\[
\lt|
\bigcup_{i=1}^{K}
P_2(R_j)\cap P_2(\wt{R}_i)
\rt|
\ge\frac13|P_2(R_j)|.
\]
Thanks to the fact that 
$|P_2(R_j)|\le|P_1(R_j)|\le|P_1(\wt{R}_i)|$, 
this implies that 
\[
\lt|\bigcup_{i=1}^{K}Q\cap\wt{R}_i\rt|
\ge\frac13|Q|,
\]
where $Q$ is a unique dyadic cube 
containing $x$ and having the side length $|P_2(R_j)|$. 
This proves \eqref{2.3}.

It follows from \eqref{2.3} and \eqref{1.2} that 
\begin{align*}
w\lt(\bigcup_{i=1}^{M}R_i\rt)
&\le 
w\lt(\lt\{
x\in\R^2:\,
\kM_{\cQ}[\1_{\bigcup_{i=1}^{N}\wt{R}_i}](x)
\ge\frac13
\rt\}\rt)
\\ &\le C
U\lt(\bigcup_{i=1}^{N}\wt{R}_i\rt)
\le C
\sum_{i=1}^{N}U(\wt{R}_i),
\end{align*}
where $U=\kM_{\cQ}w$. 
We shall evaluate the quantity 
\[
\text{(i)}
=
\sum_{i=1}^{N}U(\wt{R}_i).
\]

Let 
$\mu_{U}(x)$ be the weighted multiplicity function 
associated to the family 
$\{\wt{R}_i\}$, that is, 
\[
\mu_{U}(x)
=
\sum_{i=1}^{N}
\frac{U(\wt{R}_i)}{|\wt{R}_i|}
\1_{\wt{R}_i}(x).
\]
By \eqref{2.1}, 
choosing $\dl_0$ small enough determined later,
\begin{align*}
\text{(i)}
&\le
\sum_{i=1}^{N}
\frac{U(\wt{R}_i)}{|\wt{R}_i|}
\int_{\wt{R}_i}\frac{|f(y)|}{t}\,dy
\\ &=
\dl_0
\int_{\R^2}
\mu_{U}(x)W(x)^{-1}
\cdot
\frac{|f(x)|}{\dl_0t}
\cdot
W(x)\,dx.
\end{align*}
Using the elementary inequality
\[
ab\le(e^a-1)+b(1+\log^{+}b),
\quad a,b\ge 0,
\]
we get
\begin{align*}
\text{(i)}
&\le
\dl_0
\int_{\R^2}
\lt(\exp(\mu_{U}(x)W(x)^{-1})-1\rt)
W(x)\,dx
\\ &\quad+
\dl_0
\int_{\R^2}
\frac{|f(x)|}{\dl_0t}
\lt(1+\log^{+}\frac{|f(x)|}{\dl_0t}\rt)
W(x)\,dx
\\ &\le
\dl_0
\int_{\R^2}
\lt(\exp(\mu_{U}(x)W(x)^{-1})-1\rt)
W(x)\,dx
\\ &\quad+
(1-\log\dl_0)
\int_{\R^2}
\frac{|f(x)|}{t}
\lt(1+\log^{+}\frac{|f(x)|}{t}\rt)
W(x)\,dx.
\end{align*}

We have to evaluate the quantity 
\[
\text{(ii)}
=
\int_{\R^2}
\lt(\exp(\mu_{U}(x)W(x)^{-1})-1\rt)
W(x)\,dx.
\]
We expand the exponential in a Taylor series. 
Then 
\begin{align*}
\text{(ii)}
&=
\sum_{k=1}^{\8}
\frac1{k!}
\int_{\R^2}
(\mu_{U}(x)W(x)^{-1})^k
W(x)\,dx
\\ &=
\sum_{k=1}^{\8}
\frac1{k!}
\int_{\R^2}
\mu_{U}(x)^kW(x)^{1-k}
\,dx.
\end{align*}
Fix $k\ge 2$. 
We use an elementary inequality
\[
\lt(\sum_{i=1}^{\8}a_i\rt)^k
\le k
\sum_{i=1}^{\8}
a_i
\lt(\sum_{j=1}^ia_j\rt)^{k-1},
\]
where $\{a_i\}_{i=1}^{\8}$ 
is a sequence of summable nonnegative reals. 
Then 
\begin{align*}
\lefteqn{
\int_{\R^2}
\mu_{U}(x)^kW(x)^{1-k}
\,dx
}\\ &\le k
\sum_{i=1}^{N}
\frac{U(\wt{R}_i)}{|\wt{R}_i|}
\int_{\wt{R}_i}
\lt(
\sum_{j=1}^i
\frac{U(\wt{R}_j)}{|\wt{R}_j|}
\1_{\wt{R}_j}(x)
\rt)^{k-1}
W(x)^{1-k}\,dx
\\ &\le k
\sum_{i=1}^{N}
\frac{U(\wt{R}_i)}{|\wt{R}_i|}
\int_{\wt{R}_i}
\lt(
\sum_{j=1}^i\1_{\wt{R}_j}(x)
\rt)^{k-1}\,dx,
\end{align*}
where we have used 
\[
\sum_{j=1}^i
\frac{U(\wt{R}_j)}{|\wt{R}_j|}
\1_{\wt{R}_j}(x)
\le
\lt(\sum_{j=1}^i\1_{\wt{R}_j}(x)\rt)W(x).
\]

We claim that, 
for $n=1,2,\ldots,N$,
\begin{equation}\label{2.5}
|X_{i,n}|\le3^{1-n}|\wt{R}_i|,
\end{equation}
where 
\[
X_{i,n}
=
\lt\{
x\in\wt{R}_i:\,
\sum_{j=1}^i\1_{\wt{R}_j}(x)\ge n
\rt\}.
\]
Indeed, first we notice that, 
for any $k$ and $j$ with 
$N\ge k>j\ge 1$, if 
$\wt{R}_k\cap\wt{R}_j\ne\0$, 
then 
\[
\wt{R}_k\cap\wt{R}_j
=
P_1(\wt{R}_k)\times P_2(\wt{R}_j),
\]
because we have 
$P_1(\wt{R}_k)\subset P_1(\wt{R}_j)$ 
and, by \eqref{2.2}, 
$|P_2(\wt{R}_k)\cap P_2(\wt{R}_j)|<\frac13|P_2(\wt{R}_k)|$. 
With this in mind, we can observe the following: 

There exists a set of dyadic intervals 
$\{I_{j\,k}\}$ with 
$j=1,2,\ldots,n$ and 
$k=1,2,\ldots,K_j$ 
that satisfies the following:

\begin{itemize}
\item The dyadic intervals $I_{j\,k}$ 
are pairwise disjoint for varying $k$;
\item For each $I_{j\,k}$, $j>1$, 
there exists a unique 
$I_{(j-1)\,l}\supsetneq I_{j\,k}$;
\item For each $I_{j\,k}$ 
there exists a unique number 
$i_{j\,k}\le i$ such that 
$I_{j\,k}=P_2(\wt{R}_{i_{j\,k}})$; 
\item 
$P_2(X_{i,1})=I_{1\,1}$,
$P_2(X_{i,2})=\bigcup_{k=1}^{K_2}I_{2\,k}$,
$\ldots$,
$P_2(X_{i,j})=\bigcup_{k=1}^{K_j}I_{j\,k}$,
$\ldots$,
$P_2(X_{i,n})=\bigcup_{k=1}^{K_n}I_{n\,k}$;
\item If $I_{j\,k}\subset I_{(j-1)\,l}$, then 
$i_{j\,k}<i_{(j-1)\,l}$ and 
$\wt{R}_{i_{(j-1)\,l}}\cap\wt{R}_{i_{j\,k}}\ne\0$.
\end{itemize}

It follows from the last relation and \eqref{2.2} that 
\[
3\sum_{k=1}^{K_j}|I_{j\,k}|
<
\sum_{k=1}^{K_{j-1}}|I_{(j-1)\,k}|,
\quad j=2,3,\ldots,n.
\]
This gives us that 
\[
3^{n-1}
\sum_{k=1}^{K_n}|I_{n\,k}|
<
|I_{1\,1}|,
\]
which yields \eqref{2.5}. 

It follows from \eqref{2.5} that 
\begin{align*}
\lefteqn{
\frac{U(\wt{R}_i)}{|\wt{R}_i|}
\int_{\wt{R}_i}
\lt(
\sum_{j=1}^i\1_{\wt{R}_j}(x)
\rt)^{k-1}\,dx
}\\ &\le
\frac{U(\wt{R}_i)}{|\wt{R}_i|}
\sum_{n=1}^{N}
n^{k-1}|X_{i,n}|
\\ &\le 
U(\wt{R}_i)
\sum_{n=1}^{N}
n^{k-1}3^{1-n}.
\end{align*}
Altogether, the quantity (ii) 
can be majorized by 
\[
\text{(i)}
\times
\lt[
1+\sum_{k=2}^{\8}
\frac1{(k-1)!}
\sum_{n=1}^{N}
n^{k-1}3^{1-n}
\rt].
\]
There holds 
\[
\sum_{k=2}^{\8}
\frac1{(k-1)!}
\sum_{n=1}^{N}
n^{k-1}3^{1-n}
\le
3\sum_{n=1}^{\8}
\lt(\frac{e}{3}\rt)^n
=:C_0.
\]
If we choose $\dl_0$ so that 
$\dl_0(1+C_0)=\frac12$,
we obtain 
\[
\text{(i)}
\le C
\int_{\R^2}
\frac{|f(x)|}{t}
\lt(1+\log^{+}\frac{|f(x)|}{t}\rt)
W(x)\,dx.
\]
This completes the proof. 

\begin{xrem}
Since our proof relies only upon the dyadic structure \eqref{2.4}, 
it can be applied the basis $\cR$ of the form 
the set of all rectangles in $\R^n$ 
whose sides parallel to the coordinate axes and 
which are congruent to the rectangle 
$(0,a)^{n-1}\times(0,b)$ 
with varying $a,b>0$.
\end{xrem}

\section{Proof of Theorem \ref{thm1.4}}\label{sec3}
In what follows we shall prove Theorem \ref{thm1.4}. 
We follow the argument in \cite{G} Chapter 10, Theorem 10.3.5.
To avoid problems with antipodal points, 
it is convenient to split 
$\Sg_{N}$ as the union of eight sets, 
in each of which the angle between any two vectors 
does not exceed $\pi/4$. 
It suffices therefore to obtain 
\eqref{1.4} for each such subset of $\Sg_{N}$. 
Let us fix one such subset of $\Sg_{N}$, 
which we call $\Sg_{N}^1$. 

To prove \eqref{1.4}, 
we fix a $t>0$ and we start with 
a compact subset $K$ of the set 
$\{x\in\R^2:\,\kM_{\Sg_{N}^1}f(x)>t\}$. 
Then for every $x\in K$, 
there exists an open rectangle $R_x$ 
that contains $x$ and whose longest side
is parallel to a vector in $\Sg_{N}^1$. 
By compactness of $K$, 
there exists a finite subfamily 
$\{R_{\al}\}_{\al\in\cA}$ 
of the family $\{R_x\}_{x\in K}$ 
such that
\begin{equation}\label{3.1}
\fint_{R_{\al}}|f(y)|\,dy>t
\end{equation}
for all $\al\in\cA$ and such that 
the union of the $R_{\al}$'s covers $K$. 

In the sequel 
we denote by $\theta_{\al}$ 
the angle between the $x_1$-axis and 
the vector pointing in the longer direction of $R_{\al}$ 
for any $\al\in\cA$. 
We also denote by $l_{\al}$ 
the shorter side of $R_{\al}$ 
and by $L_{\al}$ 
the longer side of $R_{\al}$ 
for any $\al\in\cA$. 

We shall select the subfamily 
$\{R_{\bt}\}_{\bt\in\cB}$ 
as follows:

Without loss of generality 
we may assume that 
$\cA=\{1,2,\ldots,\ell\}$ with 
$L_j\ge L_{j+1}$ for all 
$j=1,2,\ldots,\ell-1$. 
Let $\bt_1=1$ and choose $\bt_2$ 
to be the first number in 
$\{\bt_1+1,\bt_1+2,\ldots,\ell\}$ 
such that 
\[
|R_{\bt_1}\cap R_{\bt_2}|
\le
\frac12|R_{\bt_2}|.
\]
We next choose $\bt_3$ 
to be the first number in 
$\{\bt_2+1,\bt_2+2,\ldots,\ell\}$ 
such that 
\[
|R_{\bt_1}\cap R_{\bt_3}|
+
|R_{\bt_2}\cap R_{\bt_3}|
\le
\frac12|R_{\bt_3}|.
\]
Suppose we have chosen the numbers 
$\bt_1,\bt_2,\ldots,\bt_{j-1}$. 
Then we choose $\bt_j$ 
to be the first number in 
$\{\bt_{j-1}+1,\bt_{j-1}+2,\ldots,\ell\}$ 
such that 
\begin{equation}\label{3.2}
\sum_{k=1}^{j-1}
|R_{\bt_k}\cap R_{\bt_j}|
\le
\frac12|R_{\bt_j}|.
\end{equation}
Since the set $\cA$ is finite, 
this selection stops after $m$ steps. 

Define 
$\cB=\{\bt_1,\bt_2,\ldots,\bt_m\}$ 
and let 
\[
Y(x)
=
\sum_{\bt\in\cB}\1_{(R_{\bt})^{*}}(x),
\]
where $(R_{\bt})^{*}$ is 
the rectangle $R_{\bt}$ expanded $5$ times 
in both directions. 

We claim that 
\begin{equation}\label{3.3}
w(K)
\le
w\lt(\bigcup_{\al\in\cA}R_{\al}\rt)
\le C(\log N)
\int_{\R^2}Y(x)U(x)\,dx,
\end{equation}
where $U(x)=\kM_{\cQ}w(x)$. 
To verify this claim, 
we need the following lemma. 

We set 
$\om_k=\frac{2\pi2^k}{N}$
for $k\in\Z^{+}$ and $\om_0=0$. 
We let 
$M=[\frac{\log(N/8)}{\log2}]$. 

\begin{lem}[\text{\cite[Lemma 10.3.6]{G}}]\label{lem3.1}
Let $R_{\al}$ be a rectangle in the family 
$\{R_{\al}\}_{\al\in\cA}$
and let $0\le k<M$.
Suppose that $\bt\in\cB$ is such that 
\[
\om_k
\le
|\theta_{\al}-\theta_{\bt}|
<
\om_{k+1}
\]
and such that 
$L_{\bt}\ge L_{\al}$. 
Let 
\[
s_{\al}
=
8\max(l_{\al},\om_kL_{\al}).
\]
For an arbitrary $x\in R_{\al}$, 
let $Q$ be a square centered at $x$ 
with sides of length $s_{\al}$ 
parallel to the sides of $R_{\al}$. 
Then we have 
\begin{equation}\label{3.4}
\frac{|R_{\bt}\cap R_{\al}|}{|R_{\al}|}
\le 32
\frac{|(R_{\bt})^{*}\cap Q|}{|Q|}.
\end{equation}
\end{lem}

We shall prove \eqref{3.3}. 
By \eqref{1.2} 
it suffices to show that 
\begin{equation}\label{3.5}
\bigcup_{\al\in\cA}R_{\al}
\subset
\lt\{x\in\R^2:\,
\kM_{\cQ}Y(x)>\frac{C}{\log N}
\rt\}.
\end{equation}
Since we may assume that 
$C/(\log N)<1$, 
the set 
$\bigcup_{\bt\in\cB}R_{\bt}$
is contained in the set of the right hand side of \eqref{3.5}. 
So, 
we fix $\al\in\cA\setminus\cB$. 
Then the rectangle $R_{\al}$ was not selected in the selection procedure. 

By the construction and \eqref{3.2}, 
we see that there exists $j$ such that 
\[
\sum_{k=1}^j
|R_{\bt_k}\cap R_{\al}|
>
\frac12|R_{\al}|
\]
and such that 
$L_{\bt_k}\ge L_{\al}$ 
for all $k=1,2,\ldots,j$. 

Let 
$\cB_j=\{\bt_1,\bt_2,\ldots,\bt_j\}$. 
It follows from Lemma \ref{lem3.1} 
that 
\begin{align*}
\frac12
&<
\sum_{\bt\in\cB_j}
\frac{|R_{\bt}\cap R_{\al}|}{|R_{\al}|}
\\ &=
\sum_{k=0}^{M}
\sum_{\substack{
\bt\in\cB_j: \\ 
\om_k\le|\theta_{\al}-\theta_{\bt}|<\om_{k+1}
}}
\frac{|R_{\bt}\cap R_{\al}|}{|R_{\al}|}
\\ &\le 32
\sum_{k=0}^{M}
\sum_{\substack{
\bt\in\cB_j: \\ 
\om_k\le|\theta_{\al}-\theta_{\bt}|<\om_{k+1}
}}
\frac{|(R_{\bt})^{*}\cap Q_k|}{|Q_k|},
\end{align*}
where $Q_k$ is a square determined by 
Lemma \ref{lem3.1} with 
an arbitrary $x\in R_{\al}$. 
Since we have 
$M\le C(\log N)$ and 
\[
\sum_{\substack{
\bt\in\cB_j: \\ 
\om_k\le|\theta_{\al}-\theta_{\bt}|<\om_{k+1}
}}
\frac{|(R_{\bt})^{*}\cap Q_k|}{|Q_k|}
\le C
\kM_{\cQ}Y(x)
\text{ for all }
x\in R_{\al},
\]
we obtain 
\[
\kM_{\cQ}Y(x)>\frac{C}{\log N}
\text{ for all }
x\in R_{\al},
\]
which implies \eqref{3.5} 
and, hence, \eqref{3.3}. 

We now evaluate 
\[
\text{(i)}
=
\int_{\R^2}Y(x)U(x)\,dx
=
\sum_{\bt\in\cB}U((R_{\bt})^{*}).
\]
By \eqref{3.1} and H\"{o}lder's inequality 
we have 
\begin{align*}
\text{(i)}
&\le
\frac{1}{t}
\sum_{\bt\in\cB}
U((R_{\bt})^{*})
\fint_{R_{\bt}}|f(y)|\,dy
\\ &=
\frac{1}{t}
\int_{\R^2}
\lt(\sum_{\bt\in\cB}
\frac{U((R_{\bt})^{*})}{|R_{\bt}|}
\1_{R_{\bt}}(y)\rt)
|f(y)|\,dy
\\ &\le
\frac{1}{t}
\lt(
\int_{\R^2}
\lt(\sum_{\bt\in\cB}
\frac{U((R_{\bt})^{*})}{|R_{\bt}|}
\1_{R_{\bt}}(y)\rt)^2
W(y)^{-1}\,dy
\rt)^{1/2}
\|f\|_{L^2(\R^2,W)}.
\end{align*}
We have further 
\begin{align*}
\text{(ii)}
&=
\int_{\R^2}
\lt(\sum_{\bt\in\cB}
\frac{U((R_{\bt})^{*})}{|R_{\bt}|}
\1_{R_{\bt}}(y)\rt)^2
W(y)^{-1}\,dy
\\ &=
\sum_{j=1}^m
\lt(\frac{U((R_{\bt_j})^{*})}{|R_{\bt_j}|}\rt)^2
\int_{R_{\bt_j}}W(y)^{-1}\,dy
\\ &\quad+2
\sum_{j=1}^m
\frac{U((R_{\bt_j})^{*})}{|R_{\bt_j}|}
\sum_{k=1}^{j-1}
\frac{U((R_{\bt_k})^{*})}{|R_{\bt_k}|}
\int_{R_{\bt_k}\cap R_{\bt_j}}
W(y)^{-1}\,dy.
\end{align*}
We notice that, for any 
$y\in R_{\bt_k}\cap R_{\bt_j}$
\[
W(y)
\ge
\frac{U((R_{\bt_k})^{*})}{|(R_{\bt_k})^{*}|}
=
\frac{U((R_{\bt_k})^{*})}{25|R_{\bt_k}|}.
\]
This yields 
\begin{align*}
\text{(ii)}
&\le 25
\sum_{j=1}^m
U((R_{\bt_j})^{*})
\\ &\quad+50
\sum_{j=1}^m
\frac{U((R_{\bt_j})^{*})}{|R_{\bt_j}|}
\sum_{k=1}^{j-1}
|R_{\bt_k}\cap R_{\bt_j}|
\\ &\le 50
\sum_{\bt\in\cB}
U((R_{\bt_j})^{*}),
\end{align*}
where we have used \eqref{3.2}. 
Altogether, we obtain 
\[
\text{(i)}
\le
\frac{C}{t^2}
\|f\|_{L^2(\R^2,W)}^2,
\]
which yields by \eqref{3.3} 
\[
w(K)
\le\frac{C(\log N)}{t^2}
\|f\|_{L^2(\R^2,W)}^2.
\]
Since $K$ was an arbitrary compact subset of 
$\{x\in\R^2:\,\kM_{\Sg_{N}^1}f(x)>t\}$,
the same estimate is valid for the latter set 
and we finish the proof.

\end{document}